\documentclass{article}
\usepackage{amssymb}

\newtheorem{theorem}{Theorem}
\newtheorem{acknowledgement}[theorem]{Acknowledgement}

\input{tcilatex}

\begin{document}

\begin{center}
\textbf{Inequalities Becker-Stark at extreme points}\bigskip

Ling Zhu$^{1}$ and Cristinel Mortici$^{2}$\bigskip

$^{1}$Department of Mathematics, Zhejiang Gongshang University, Hangzhou
310018, China, zhuling0571@163.com

$^{2}$Valahia University of T\^{a}rgovi\c{s}te, Bd. Unirii 18, 130082, T\^{a}%
rgovi\c{s}te, Romania, cristinel.mortici@hotmail.com

\[
\]
\end{center}

\textbf{Abstract: }\emph{The aim of this work is to extend Becker-Stark
inequalities near the origin and }$\pi /2.$\bigskip



\section{Introduction and Motivation}

In 1978, M. Becker and L. E. Stark presented the inequalities%
\begin{equation}
\frac{8}{\pi ^{2}-4x^{2}}<\frac{\tan x}{x}<\frac{\pi ^{2}}{\pi ^{2}-4x^{2}}%
,\ \ \ 0<x<\frac{\pi }{2}.  \label{a}
\end{equation}%
These inequalities were proven to be of great interest through the
researchers, since they were extended in different forms in the recent past.
We refer to \cite{bk}-\cite{zzz} and all references therein.

It is true that inequalities (1) and some of recent improvements are nice
through its symmetric form, but let us remember the practical importance of
an inequality which is to provide some bounds for a given expression.
Observe that near $\pi /2,$ the right-hand side inequality (1) becomes weak,
as soon as%
\[
\lim_{x\rightarrow \left( \pi /2\right) _{-}}\left( \frac{\pi ^{2}}{\pi
^{2}-4x^{2}}-\frac{\tan x}{x}\right) =\infty . 
\]%
The left-hand side inequality (1) is somehow motivated by the fact that%
\[
\lim_{x\rightarrow \left( \pi /2\right) _{-}}\left( \pi ^{2}-4x^{2}\right) 
\frac{\tan x}{x}=8. 
\]%
Having in mind these remarks, we assume that good estimates of $\left( \tan
x\right) /x$ near $\pi /2$ are obtained using expressions of the form%
\[
\frac{\tan x}{x}\approx \frac{8+\omega \left( x\right) }{\pi ^{2}-4x^{2}}\
,\ \ \ x<\frac{\pi }{2}, 
\]%
with $\omega \left( x\right) $ tending to zero, as $x$ approaches $\pi /2,$
with $x<\pi /2.$

More precisely, we propose the following double inequality of Becker-Stark
type, on a neighborhood of $\pi /2.$\bigskip

\textbf{Theorem 1. }\emph{For every }$x\in \left( 0.373,\pi/2\right) $ \emph{%
in the left-hand side and for every }$x\in \left( 0.301,\pi/2\right) $ \emph{%
in the right-hand side}$,$ \emph{the following inequalities hold true:}%
\begin{equation}
\frac{8+a\left( x\right) }{\pi ^{2}-4x^{2}}<\frac{\tan x}{x}<\frac{8+b\left(
x\right) }{\pi ^{2}-4x^{2}},  \label{t}
\end{equation}%
\emph{where}%
\[
a\left( x\right) =\frac{8}{\pi }\left( \frac{\pi }{2}-x\right) +\left( \frac{%
16}{\pi ^{2}}-\frac{8}{3}\right) \left( \frac{\pi }{2}-x\right) ^{2} 
\]%
\emph{and}%
\[
b\left( x\right) =a\left( x\right) +\left( \frac{32}{\pi ^{3}}-\frac{8}{3\pi 
}\right) \left( \frac{\pi }{2}-x\right) ^{3}. 
\]

We think to similar comments on the behavior of left-hand inequality (1)
near the origin. In connection with this problem, we propose the following
improvement.\bigskip

\textbf{Theorem 2. }\emph{For every real number }$x\in \left( 0,1.371\right)
,$ \emph{the following inequality holds true:}%
\begin{equation}
\frac{\tan x}{x}<\frac{\pi ^{2}-\left( 4-\frac{1}{3}\pi ^{2}\right)
x^{2}-\left( \frac{4}{3}-\frac{2}{15}\pi ^{2}\right) x^{4}}{\pi ^{2}-4x^{2}}.
\label{r}
\end{equation}

As a fact to remember, if we wish to obtain accurate approximations of $%
\left( \tan x\right) /x,$ using Becker-Stark inequalities, then the best
constant at the numerator near $0$ is $\pi ^{2},$ while the best constant at
the numerator near $\pi /2$ is $8.$

\section{The Proofs}

\emph{Proof of Theorem 1. }Inequalities (\ref{t}) are equivalent to $f<0$ on 
$\left( 0.373,\pi/2 \right) $ and $g>0$ on $\left( 0.301,\pi/2 \right) ,$
where%
\[
f\left( x\right) =\arctan \left( x\cdot \frac{8+a\left( x\right) }{\pi
^{2}-4x^{2}}\right) -x\text{ \ \ and \ \ \ }g\left( x\right) =\arctan \left(
x\cdot \frac{8+b\left( x\right) }{\pi ^{2}-4x^{2}}\right) -x. 
\]%
If $p_{k},$ $q_{k}$ are polynomial functions, we use the following
derivation formula%
\begin{equation}
\left( \arctan \frac{p_{k}\left( x\right) }{q_{k}\left( x\right) }-x\right)
^{\prime }=\frac{p_{k}{}^{\prime }\left( x\right) q_{k}\left( x\right)
-p_{k}\left( x\right) q_{k}{}^{\prime }\left( x\right) -p_{k}^{2}\left(
x\right) -q_{k}^{2}\left( x\right) }{p_{k}^{2}\left( x\right)
+q_{k}^{2}\left( x\right) }.  \label{d}
\end{equation}%
With $p_{1}\left( x\right) =x\left( 8+a\left( x\right) \right) $ and $%
q_{1}\left( x\right) =\pi ^{2}-4x^{2},$ we get%
\begin{equation}
f^{\prime }\left( x\right) =\frac{\left( \pi -2x\right) ^{3}u\left( x\right) 
}{9\pi ^{4}\left( p_{1}^{2}\left( x\right) +q_{1}^{2}\left( x\right) \right) 
},  \label{f}
\end{equation}%
where%
\begin{eqnarray*}
u\left( x\right) &=&\allowbreak 144\pi ^{3}-15\pi ^{5}+x\left( 432\pi
^{2}-42\pi ^{4}\right) \\
&&+\allowbreak x^{2}\left( 96\pi ^{3}-432\pi -4\pi ^{5}\right) +x^{3}\left(
8\pi ^{4}-96\pi ^{2}+288\right) .
\end{eqnarray*}%
Now we present the following steps in our proof:

\begin{itemize}
\item $u^{\prime \prime }$ is a first degree polynomial funcion, strictly
increasing, with $u^{\prime \prime }\left( 0.373\right) =\allowbreak
1058.\,\allowbreak 803...>0,$ so $u^{\prime \prime }>0$ on $\left(
0.373,\pi/2\right) .$

\item $u^{\prime }$ is strictly increasing, with $u^{\prime }\left(
0.373\right) =\allowbreak 517.\,\allowbreak 421...>0,$ so $u^{\prime }>0$ on 
$\left( 0.373,\pi/2 \right) .$

\item $u$ is strictly increasing, with $u\left( 0.373\right) =\allowbreak
0.168...,$ so $u>0$ on $\left( 0.373,\pi/2\right) .$
\end{itemize}

By (\ref{f}), $f$ is strictly increasing, with limit $f\left(\pi/2-0\right)
=0,$ so $f<0$ and the first inequality (\ref{t}) is proved.

With $p_{2}\left( x\right) =x\left( 8+b\left( x\right) \right) $ and $%
q_{2}\left( x\right) =\pi ^{2}-4x^{2}$ in (\ref{d})$,$ we get%
\begin{equation}
g^{\prime }\left( x\right) =-\frac{\left( \pi -2x\right) ^{4}v\left(
x\right) }{9\pi ^{6}\left( p_{2}^{2}\left( x\right) +q_{2}^{2}\left(
x\right) \right) },  \label{g}
\end{equation}%
where%
\begin{eqnarray*}
v\left( x\right) &=&\allowbreak 18\pi ^{6}-180\pi ^{4}+x\left( 60\pi
^{5}-576\pi ^{3}\right) +x^{4}\left( 4\pi ^{4}-96\pi ^{2}+576\right) \\
&&+\allowbreak x^{3}\left( 240\pi ^{3}-1152\pi -12\pi ^{5}\right)
+x^{2}\left( 864\pi ^{2}-168\pi ^{4}+9\pi ^{6}\right) \allowbreak .
\end{eqnarray*}%
Now we present the following steps in our proof:

\begin{itemize}
\item $v^{\prime \prime }$ is a second degree polynomial with minimum at $%
x_{0}=-2.\,\allowbreak 067....$ As $v^{\prime \prime }\left( 0.301\right)
=1921.\,\allowbreak 145...>0,$ it results that $v^{\prime \prime }>0$ on $%
\left( 0.301,\pi/2\right) .$

\item $v^{\prime }$ is strictly increasing, with $v^{\prime }\left(
0.301\right) =\allowbreak 1035.\,\allowbreak 057...>0,$ so $v^{\prime }>0$
on $\left( 0.301,\pi/2 \right) .$

\item $v$ is strictly increasing, with $v\left( 0.301\right) =\allowbreak
0.434\,386\,67>0,$ so $v>0$ on $\left( 0.301,\pi/2 \right) .$
\end{itemize}

By (\ref{g}), $g$ is strictly decreasing, with limit $g\left( \pi/2-0\right)
=0,$ so $g>0$ on $\left( 0.301,\pi/2 \right) $\ and Theorem 1 is proved.$%
\square $\bigskip

\emph{Proof of Theorem 2. }We rewrite (\ref{r}) in the form $h>0$ on $\left(
0,1.371\right) ,$ where%
\[
h\left( x\right) =\arctan \left( x\cdot \frac{\pi ^{2}-\left( 4-\frac{1}{3}%
\pi ^{2}\right) x^{2}-\left( \frac{4}{3}-\frac{2}{15}\pi ^{2}\right) x^{4}}{%
\pi ^{2}-4x^{2}}\right) -x. 
\]%
With $p_{3}\left( x\right) =x\left( \pi ^{2}-\left( 4-\frac{1}{3}\pi
^{2}\right) x^{2}-\left( \frac{4}{3}-\frac{2}{15}\pi ^{2}\right)
x^{4}\right) $ and $q_{3}\left( x\right) =\pi ^{2}-4x^{2}$ in (\ref{d}), we
get

\begin{equation}
h^{\prime }\left( x\right) =-\frac{x^{6}w\left( x^{2}\right) }{225\left(
p_{3}^{2}\left( x\right) +q_{3}^{2}\left( x\right) \right) },  \label{az}
\end{equation}%
where

\[
w\left( t\right) =85\pi ^{4}-840\pi ^{2}+t^{2}\left( 4\pi ^{4}-80\pi
^{2}+400\right) +t\left( 20\pi ^{4}-440\pi ^{2}+2400\right) . 
\]%
This function $w$ is a second degree polynomial function with minimum at $%
t_{0}=-40.\,\allowbreak 844...$ . As $w\left( 1.881\right) =\allowbreak
-0.0037...<0,$ it follows that $w\left( t\right) <0$ for every $t\in \left(
0,1.881\right) .$ That is $w\left( x^{2}\right) <0,$ for every $x\in \left(
0,1.371\right) $ (Notice that $\sqrt{1.881}=\allowbreak 1.\,\allowbreak
371\,495...$).

By (\ref{az}), $h$ is strictly increasing on $\left( 0,1.371\right) ,$ with $%
h\left( 0\right) =0,$ so $h\left( x\right) >0,$ for every $x\in \left(
0,1.371\right) .$%
\[
\]

\begin{acknowledgement}
The work of the second author was supported by a grant of the Romanian
National Authority for Scientific Research, CNCS-UEFISCDI project number
PN-II-ID-PCE-2011-3-0087.
\end{acknowledgement}

\end{document}